\documentclass{amsart}
\usepackage[pagebackref]{hyperref}

\newtheorem{lemma}{Lemma}

\newtheorem{theorem}[lemma]{Theorem}
\newtheorem{cor}[lemma]{Corollary}

\newcommand{\R}{{\mathbb R}}
\newcommand{\Z}{{\mathbb Z}}
\newcommand{\Q}{{\mathbb Q}}
\newcommand{\Log}{\mbox{Log}\,}

\newcommand{\eps}{\varepsilon}
\newcommand{\ph}{\phantom{-}}
\newcommand{\OO}{{\mathcal O}}

\title{Simplest Cubic Fields}
\author{Franz Lemmermeyer, Attila Peth\"o}

\begin{document}

\maketitle

\section{Introduction}

In this paper we intend to show that certain integers do not 
occur as the norms of principal ideals in a family of 
cubic fields studied by Cohn \cite{Cohn}, Shanks \cite{Shanks},
and Ennola \cite{Ennola}. 
These results will simplify the construction of 
certain unramified quadratic extensions of such fields (cf. 
\cite{Washington}, \cite{Watabe} etc.). 

For a natural number $a$, let $f_{a}=x^{3}-ax^{2}-(a+3)x-1$, and
let $K=K_{a}$ be the cyclic cubic number field generated by a root
$\alpha$ of $f_{a}$. Let N denote the Norm $N_{K/\Q}$.
Elements in
$K$ are said to be {\em associated} if their quotient is a unit in
$\Z[\alpha]$. The polynomial $f=f_{a}$ has discriminant
disc $ f = m^{2}$, where $m = a^{2}+3a+9$; if we assume $m$ to be
squarefree, then we have disc $\,K\,=$ disc $ f = m^{2}$ and
$\OO_{K} = \Z[\alpha]$ (there exist infinitely many such $m$, cf. 
Cusick \cite{Cusick}). Moreover it is easy to see that
$\{1 \, ,\alpha \, ,\alpha'\}$ also is an integral basis of
$\OO_{K}$: in fact, this follows from
$(\alpha+1)(\alpha^{2}-(a+1)\alpha-2)=-1$. For this family of
cyclic cubic fields, we will prove the following result:
\begin{theorem}\label{th6}
For all $\gamma \in \Z[\alpha]$ either $|N\gamma| \ge 2a+3$,
or $\gamma$ is associated to an integer. Moreover, if $|N\gamma|=2a+3$, 
then $\gamma$ is associated to one of the conjugates of $\alpha-1$.
\end{theorem}

\section{The Proof}

We start the proof with the observation that the assertion is correct
for $a<7$ ("proof by inspection" using the decomposition law for cyclic 
cubic fields or by using the method described below, but with the actual 
values of $\alpha, \alpha'$ and 
$\alpha''$). Moreover, we remark that
$\alpha$, $\alpha'=-\frac{\alpha+1}{\alpha}$
and $\alpha''=-\frac{1}{\alpha+1}$
are the roots of $f$. Choosing $\alpha$ as the smallest of the three
roots and applying Newton's method, we find that
\[-1-\frac{1}{a} < \alpha < -1-\frac{1}{2a},\quad
  -\frac{1}{a+2} < \alpha'< -\frac{1}{a+3},\quad
  a+1 < \alpha'' < a+1+\frac{2}{a}.\]
These inequalities imply (for $a \ge 7$)
\[|\alpha-\alpha'| < 1+\frac{1}{a}, \quad
  |\alpha'-\alpha''|<a+1+\frac{3}{a},\quad
\mbox{ and } |\alpha''-\alpha| < a+2+\frac{3}{a}.\]
In particular, we have
\[|\alpha-\alpha'|+
  |\alpha'-\alpha''|+
  |\alpha''-\alpha| < 2a+4+\frac{7}{a} \le 2a+5.\]
Moreover, we will need the relation
\begin{displaymath}
m=\alpha^{2}+\alpha'^{2}+\alpha''^{2}-\alpha\alpha'-\alpha'\alpha''
  -\alpha''\alpha = \mbox{det } \left( \begin{array}{ccc}
     1 & 1 & 1 \\
     \alpha & \alpha' & \alpha'' \\
     \alpha' & \alpha'' & \alpha
     \end{array} \right),
\end{displaymath}
which can be deduced easily from the well known fact that the
square of this determinant equals disc $(1, \alpha, \alpha')=m^{2}$,
making use of the formulae
\[\alpha\alpha'=-\alpha-1, \quad \alpha^{2}=a+2+a\alpha-\alpha'.\]
The units $\alpha'^{-1}$ and $\alpha''$ satisfy the inequalities
\[a+2 < |\alpha'^{-1}| < a+3, \quad a+1 < |\alpha''| < a+2,\]
and this implies that, given two positive real numbers $c_{1}$
and $c_{2}$ and an element $\gamma \in Z[\alpha]$, we can find
a unit $\eta$ such that
\begin{equation}\label{*}
c_{1} \le |\gamma \eta | < (a+3)c_{1}, \quad
  c_{2} \le |\gamma'\eta'| < (a+4)c_{2}.
\end{equation}

This is a special case of a more general result which is valid
for all number fields with unit rank $\ge 1$; we will, however,
give the proof only for totally real cubic fields K because the
notation simplifies considerably. Let $u_{1}$ and $u_{2}$ be two
independent units in $K$; their images in $\R^{2}$ upon their
logarithmic embedding are
    $\Log(u_{1}) = v_{1} = (\log |u_{1}|,\log |u_{1}'|)$
and $\Log(u_{2}) = v_{2} = (\log |u_{2}|,\log |u_{2}'|)$.
Dirichlet's unit theory shows that $v_{1}$ and $v_{2}$ are linear
independent vectors. This implies that, for any $\xi \in K$, its
image $\Log(\xi)$ can be moved into the fundamental domain spanned
by $v_{1}$ and $v_{2}$ by adding and subtracting suitable multiples
of $v_{1}$ and $v_{2}$, i.e. we can find a translate $\eta$ of
$Log(\xi)$ such that
$$ \begin{array}{rclrcl} \\ \medskip 
 c_{1} & < & |\eta|  & \le & c_{1} +
  \Big|\log |u_{1}|\Big| + \Big|\log |u_{2}|\Big|, \\ \medskip 
  c_{2} & < & |\eta'| & \le & c_{2} +
  \Big|\log |u_{1}'|\Big|+ \Big|\log |u_{2}'|\Big| \end{array} $$ 

Translating this back to the field $K$ and using the units $\alpha$
and $\alpha''$, we see that we can find  a unit $\eta$ such that
\[c_{1} \le |\gamma \eta | < |\alpha\alpha''|c_{1}, \quad
  c_{2} \le |\gamma'\eta'| < |\alpha'^{-1}\alpha|c_{2},\]
where we have chosen the exponents of $\alpha,\,\alpha',\,\alpha''$
in such a way that their absolute value is $> 1$ (this comes
from the absolute values on the $\log$'s). Inserting the bounds
on $|\alpha|,\,|\alpha'|,\,|\alpha''|$ we get equation (\ref{*}).

Writing $\xi=\gamma\eta=r+s\alpha + t\alpha'$ and $n=|N_{K/\Q}\xi|$,
we find ($T = T_{K/\Q}$ denotes the trace):
\[mt = T\bigl( \xi(\alpha' -\alpha'')\bigr), \quad
  ms = T\bigl( \xi(\alpha''-\alpha  )\bigr), \quad
  mr = T\bigl( \xi(\alpha\alpha'-\alpha''^{2})\bigr). \]
Letting $c_{1} = c_{2} = \sqrt[3\,]{n/(a+3)}$ we get
$|\xi|,\,|\xi'|,\,|\xi''|\,< \sqrt[3\,]{n}\cdot (a+3)^{2/3},$
and this implies the bounds
\[\begin{array}{lcl} \smallskip
|mt|  & \le &  |\xi| |\alpha'-\alpha''|+|\xi'| |\alpha''-\alpha|
              +|\xi''| |\alpha-\alpha'| \\ \smallskip
      & <   & \sqrt[3\,]{n}(a+3)^{2/3}(2a+5); \\ \smallskip
|ms|  & \le &  |\xi| |\alpha''-\alpha|+|\xi'| |\alpha-\alpha'|
              + |\xi''| |\alpha'-\alpha''|\\ \smallskip
      & <   & \sqrt[3\,]{n}(a+3)^{2/3}(2a+5).
\end{array}\]
Using $n \le 2a+3 $ and $a \ge 7$ we find that
$|t| \le 2,|s| \le 2.$ Computing the actual values of $\alpha,
\alpha'$ and $\alpha''$ for $1 \le a \le 6$ and carrying out the
above procedure we get the same result.

Now we will look at the $\xi = r+s\alpha+t\alpha'$ that satisfy
the following system of inequalities:
\[|s| \le 2, \ |t| \le 2,\quad |\xi \xi' \xi''| \le n \le 2a+3,\]
\[|\xi|, \ |\xi'|, \ |\xi''| < \sqrt[3\,]{n}(a+3)^{2/3}. \]
A somewhat tedious computation yields
$N_{K/\Q}(r+s\alpha+t\alpha')=r^{3}+s^{3}+t^{3}+ar^{2}s+ar^{2}t
 +3st^{2}-(a^{2}+3a+6)s^{2}t-(a+3)rt^{2}-(a+3)rs^{2}+(a^{2}+a+3)rst,$ 
 so for fixed $s, t$ the norm of $r+s\alpha+t\alpha'$ is a cubic polynomial
in $r$. This polynomial will be minimal for values of $r$ in the
neighborhood of its roots. We will distinguish the following cases:
\begin{enumerate}
\item $s=t=0$:
      then $\xi \in \Z$, and $\xi$ (as well as $\gamma$) is
      associated to a natural number;
\item $s=\pm 1$, $t=0$:  then $\xi = \alpha-r$ for some
      $r \in \Z$, and we find\\
      $N_{K/\Q}\xi = -f(r) = -r^{3}+ar^{2}+(a+3)r+1.$\\
      The roots of this polynomial are $r\approx 0, r\approx -1,$
      and $r\approx a+1$, and now \\
      $N\alpha = -N(\alpha+1) = 1,\\
       N(\alpha-1)= N(\alpha+2) = 2a+3, \quad
       N(\alpha-a-1)= 2a+3,\\ %
       N(\alpha-a) = -N(\alpha-a-2) = a^{2}+3a+1 > 2a+3$, if $a \ge 2$\\
       show that either $\xi$ is associated to $1$, or
       $|N\xi| \ge 2a+3$.
\item $s=\pm 2$, $t=0$: proceeding as in case 2 and keeping in
       mind that we need examine only those $\xi=2\alpha-r$ with
       $r$ odd, we find that\\
       $N(2\alpha+1)=-(2a+3), \quad N(2\alpha-1)=6a+19,$ \\
       $N(2\alpha+3)=6a-1, \quad N(2\alpha-2a-1)=4a^{2}+24a+19,$\\
       $N(2\alpha-2a-3)=-4a^{2}+17.$
\item $s=0, t = \pm 1$ and $s=0, t = \pm 2$: this yields nothing
       new, because $\xi' = r+s\alpha'$ has already been examined.
\item $s=t=\pm 1$: then $\xi = \alpha + \alpha' - r = -\alpha''
       +a -r$, and therefore $\xi' = -\alpha +a-r$ is of the
       type studied in 2; since $N\xi = N\xi'$ we are done.
\item $s=-t=\pm 1$: then $\xi = r+\alpha-\alpha'$, and
      \[f(r) = N(\xi) = r^3 - (a^2+3a+9) r +(a^2+3a+9),\]
 \[ \begin{array}{lcrlcr}
    f(1)    & = & 1,     \quad    & f(2)    & = & -a^2-3a-1,\\
    f(a+1)  & = & -6a+1, \quad    & f(a+2)  & = & 2a^2-1,  \\
    f(-a-2) & = & 6a+19, \quad    & f(-a-3) & = & -2a^2-6a+9
   \end{array} \]
\item $s=\pm 2, t=\mp 1$: then $\xi = r+2\alpha-\alpha'$, and
     \[f(r) = N(\xi) = r^3 + a r^2 - (2a^2+7a+21) r +(4a^2+12a+37), \]
\[ \begin{array}{lcrlcr}
    f(2)     & = & 2a+3,    \quad    & f(3)     & = & -2a^2+1,    \\
    f(a+1)   & = & -12a+17, \quad    & f(a+2)   & = & 3a^2-7a+3, \\
    f(-2a-3) & = & 30a+37,  \quad    & f(-2a-4) & = & -6a^2+2a+57
 \end{array} \]
\item $s=\pm 2, t=\pm 1$: then $\xi = r+2\alpha+\alpha'=r+a+\alpha-\alpha''$,
and we can proceed as in 5, refering to case 6 instead of 2.
\item $s = -t = \pm 2$: then $\xi = r+2\alpha-2\alpha'$, and according to 6 
we have to consider only the case $r$ odd, thus
      \[f(r) = N(\xi) = r^3- (12a + 36 + 4a^2)r + 8 a^2 + 24 a +72, \]
\[ \begin{array}{llll}
 f(1)      & =   4a^2+12a+37,   & f(3)      & =  -4a^2-12a-9,  \\
 f(2a+1)   & =  -8a^2-54a+37,   & f(2a+3)   & =  8a^2-30a-9,  \\
 f(-2a-3)  & =   8a^2+78a+153,  & f(-2a-5)  & =  -8a^2+6a+127.
 \end{array} \]
\item $s=t= \pm 2:$ then $\xi = r+2\alpha+2\alpha'=r+2(-\alpha''+a)=
r+2a-2\alpha''$ and we can proceed as in 5.
\end{enumerate}
 
Assume now that $|N\gamma|=2a+3$. Then the proof of the first assertion 
shows that $\gamma$ is associated to one of the elements given in the 
second column Table 1.

\begin{table}[h]
\caption{}
\vspace{0.5ex}
\hfil \begin{tabular}{|c|l|l|} \hline
     & $ \alpha-1                   $ & $                                 $ \\
     & $ \alpha+2                   $ & $ -(\alpha-1)''(\alpha+1)         $ \\
     & $ \alpha-(a+1)               $ & $ -(\alpha-1)'/(\alpha+1)         $ \\
     & $ 2 \alpha+1                 $ & $ -(\alpha-1)' \alpha  	          $ \\
     & $ \alpha +\alpha'-a+1        $ & $ -(\alpha-1)''          	  $ \\
     & $ \alpha +\alpha'-a-2        $ & $ \ph(\alpha-1)'\alpha/(\alpha+1) $ \\
     & $ \alpha +\alpha'+1          $ & $ \ph(\alpha-1)(\alpha+1)/\alpha  $ \\
     & $ 2\alpha-\alpha'+2          $ & $ -(\alpha-1)'(\alpha+1)          $ \\
     & $ 2\alpha+2\alpha'-2a-1      $ & $ -(\alpha-1)/(\alpha+1)          $ \\ 
     \hline 
$a=1 $ & $ \alpha-3              $ & $ \ph(\alpha-1)''/(\alpha+1)         $ \\
     & $ 2\alpha+3             $ & $ \ph(\alpha-1)(\alpha+1)^{2}/\alpha   $ \\
     & $ \alpha+\alpha'+2      $ & $ -(\alpha-1)'(\alpha+1)/\alpha        $ \\
     & $ \alpha-\alpha'+2      $ & $ \ph(\alpha-1)(\alpha+1)              $ \\
     & $ 2\alpha+\alpha'-3     $ & $ -(\alpha-1)''\alpha /(\alpha+1)      $ \\
     & $ 2\alpha+2\alpha'-5    $ & $  \ph(\alpha-1)''\alpha^2/(\alpha+1)  $ \\ 
     \hline 
$a=2 $ & $ \alpha-\alpha'+4     $ & $ \ph(\alpha-1)\alpha                 $ \\
     & $ 2\alpha-\alpha'+3    $ & $ \ph(\alpha-1)(\alpha+1)               $ \\
     & $ 2\alpha+\alpha'-6    $ & $ \ph(\alpha-1)''/(\alpha+1)            $ \\ 
     \hline 
$a=3 $ & $ 2\alpha-\alpha'+5    $ & $ \ph(\alpha-1)\alpha                 $ \\
     & $ 2\alpha-\alpha'-10   $ & $ -(\alpha-1)''/\alpha(\alpha+1)$ \\ \hline
\end{tabular}  \hfil \end{table}

In the third column  we have given the 
factorization  of the corresponding element as 
$\alpha-1$ or $(\alpha-1)'$ or $(\alpha-1)''$ times a unit. The 
table consists of four subtables: in the first we collected those 
numbers whose norms are $2a+3$ for any $a$. In the remaining 
subtables you find the exceptional elements which appear only for 
$a=1, 2$ or $3$. The subtables are indicated with the values of $a$. 
The proof of the identities is straightforward and therefore omitted. 

This completes the proof of Theorem \ref{th6}. We acknowledge
the help of Maple V (version 4.4) and PARI (version 1.38.3)
in checking the computations. 

\section{Applications}

From Theorem \ref{th6} we deduce the following

\begin{cor}\label{th7}
Assume that $m$ is squarefree and that $2a+3=b^{2}$ for some
$b \in \Z$. Then, $L=K(\sqrt{\alpha+2},\, \sqrt{\alpha'+2})$
is a quartic unramified extension of $K$ with
$Gal(L/K) \cong C_{2} \times C_{2}$. In particular, $Cl(L)$
contains a subgroup of type $C_{2} \times C_{2}.$
\end{cor}

\medskip\noindent{\it Proof.}
Suppose that $\alpha+2$ is a square in $\OO_K$; since
$N(\alpha+2) = 2a+3 = b^2$, this implies that there is an element
$\gamma$ of norm $b < 2a+3$. Theorem 3.1. implies that $\gamma$ is
associated to an integer $r \in \Z$, hence  $\alpha + 2 = r^2\eps$
for some unit $\eps \in \OO_K^\times$. But $\{1, \alpha,\alpha'\}$ is
an integral basis of $\OO_K$, hence $r \mid (\alpha+2)$ implies
that $r = \pm 1$.

The rest of the proof is the same as in \cite{Washington},
\cite{Watabe} or \cite{Lan}.

\medskip

Similarly, we can show (cf. \cite{Washington}):

\begin{cor}\label{th8}
Assume that $m>13$ is squarefree and that
$6a+19=b^{2}$ for some $b \in \Z$. Then,
$L = K(\sqrt{\alpha(2\alpha-1)},\,\sqrt{\alpha'(2\alpha'-1)})$
is a quartic unramified extension of $K$ with
$Gal(L/K) \cong C_{2} \times C_{2}$. In particular, $Cl(L)$
contains a subgroup of type $C_{2} \times C_{2}$.
\end{cor}

\noindent
{\bf Remark:} If $a=1$, $m=13$, we have $6a+19 = (2a+3)^{2}$, i.e.
$b$ is a norm.


\begin{thebibliography}{ACH}
\bibitem[C]{Cohn}
H. Cohn, {\em A device for generating fields of even class number},
Proc. Amer. Math. Soc. {\bf 7} (1956), 595--598
%

\bibitem[Cu]{Cusick}
T.W. Cusick, {\em Lower bounds for regulators}, in.: Number Theory
Noordwijkerhout, 1983, pp. 63--73, Ed.: H. Jager, LNM Vol. 1068,
Springer Verlag, 1984
%

\bibitem[E]{Ennola}
V. Ennola, {\em Cubic fields with exceptional units}, Computational Number 
Theory Debrecen, Hungary 1989, pp. 103--128, Eds.: A. Peth\H{o}, M.E. Pohst,
H.C. Williams and H.G. Zimmer, Walter deGruyter Verlag, 1991.
%

\bibitem[L]{Lan}
Y.-Z. Lan, {\em Arithmetic properties of a class of cyclic cubic
fields}, Sci. China {\bf 32} (1989), 922--928
%

\bibitem[Sh]{Shanks}
D. Shanks, {\em The simplest cubic number fields}, Math. Comp.
{\bf 28} (1974), 1137--1152
%

\bibitem[Wa]{Washington}
L. Washington, {\em Class numbers of the simplest cubic fields},
Math. Comp. {\bf 48} (1987), 371--384
%

\bibitem[W]{Watabe}
M. Watabe, {\em On certain cubic fields IV}, Proc. Japan Acad. {\bf 59A}
(1983), 387--389
%

\end{thebibliography}
\end{document}